\documentclass[12pt]{article}
\usepackage{graphicx}
\usepackage{psfrag}
\usepackage{amsmath}
\usepackage{amssymb}
\usepackage{amsthm}

\author{Evelyn Sander \\
	Department of Mathematical Sciences, 
	George Mason University and \\
Jim Yorke, University of Maryland
}

\date{Preprint \today}

\newtheorem{lemma}{Lemma}
\newtheorem{definition}{Definition}

\newtheorem{theorem}{Theorem}
\newtheorem{remark}{Remark}

\begin{document}

\title{A classification of explosions in dimension one}
\author{E.~Sander - J.A. Yorke
}

\maketitle
\begin{abstract}
  A discontinuous change in the size of an attractor is the most
  easily observed type of global bifurcation. More generally, an {\it
    explosion} is a discontinuous change in the set of recurrent
  points. An explosion often results from heteroclinic and homoclinic
  tangency bifurcations. Newhouse and Palis conjectured in 1976 that
  planar explosions are generically the result of either tangency
  or saddle node bifurcations. In this paper, we prove
  this conjecture for one-dimensional maps. Furthermore, we give a
  full classification for all possible tangency bifurcations and
  whether they lead to explosions.
\end{abstract}

\section{Introduction}

For continuously varying one-parameter families of iterated maps in
$R^n$, discontinuous changes in the size of an attractor are the most
easily observed type of global bifurcations, including changes in the
basin boundary (metamorphasis).  These changes can occur as the result
of a change in stability of the recurrent set. A more general
situation occurs when discontinuous changes in attractors occur as the
result of a discontinuous change in the size of the recurrent set
itself. Such a global bifurcation is called an {\bf explosion}.  For
the last several years, we have been studying explosions and their
properties, including a classification of explosions at heteroclinic
tangencies for planar diffeomorphisms~\cite{alligood:sander:yorke:02},
and more recently a numerical study of the statistical properties of a
certain kind of explosions that occur in dimension three and
higher~\cite{alligood:sander:yorke:06}.  Our research, as well as that
of many others~\cite{horita:muniz:sabini:07, sabini:01}, has been
guided by a 1976 conjecture of Newhouse and
Palis~\cite{newhouse:palis:76}. (See also the restatements
in~\cite{bonatti:diaz:viana:05, palis:takens:93a}). For over thirty
years, this conjecture has managed to elude proof. The conjecture says
that for for generic planar diffeomorphisms, all explosions occur
through the following two local bifurcations: saddle-node bifurcations
and homoclinic bifurcations.  In this paper, we prove that the
conjecture is true for smooth interval maps. In addition, we are able
to give a full classification of which types of tangencies give rise
to explosions.  These are important and new result in themselves. They
are for example in contrast with a recent result of Horita,
Muniz, and Sabini~\cite{horita:muniz:sabini:07}, showing that in a
probabilistic sense, the Newhouse Palis conjecture is not true for
circle maps. In a broader sense, we are hopeful that the insight
gained from the one-dimensional case will give rise to insights
leading to the proof of the planar case of the Newhouse Palis
conjecture. 

It is clear that an isolated saddle node bifurcation for either a
fixed point or periodic orbit gives rise to a local explosion, since
new periodic points appear. However, in many cases in both one and
higher dimensions, a saddle node bifurcation also gives rise to a
global bifurcation. The set of recurrent points changes
discontinuously as the parameter is varied at points not contained in
the saddle node periodic orbit.  For example, the period three window
in the chaotic attractor for the logistic map corresponds to an
explosion: For parameter values less than a bifurcation value, there
is a global attractor which comprises the full recurrent set. It
consists of an interval. After the bifurcation parameter, the global
attractor consists only of a period three orbit, and the full
recurrent set is a nowhere dense Cantor set within the attractor
interval. An explosion occurs at the bifurcation parameter, and it is
due to a saddle-node bifurcation at a point inside the global
attractor. Saddle node bifurcations an invariant circles are
another well-studied examples of this phenomenon.

In all dimensions, explosions due to a homoclinic or heteroclinic
bifurcation essentially occur due to the creation of a homoclinic
tangle when stable and unstable manifolds interesect transversally. In
one and two dimensions, the transition between no intersection and
transverse intersections involve tangencies between the stable and
unstable manifolds of fixed or periodic points ({\it cf.}
\cite{bonatti:diaz:viana:05}). For example, the H\'enon map and Ikeda
map contain well studied examples of explosions which are a result of
homoclinic bifurcations. In higher dimensions, such bifurcations can
occur without tangencies~\cite{diaz:rocha:02}. Such a bifurcation
leads to unstable dimension variability~\cite{alligood:sander:yorke:06}.

In our previous work for planar maps, we gave a precise
classification for which types of tangencies for heteroclinic cycles
will result in explosions. We called this class of cycles {\bf
  crossing cycles}, because the different stable and unstable
manifolds involved in the cycle lie {\it across} the tangency from
each other. We show here that the same results hold for interval
maps. Our main results are as follows:

\begin{theorem}[Explosions at tangencies]\label{T:explosions}
  For generic one-parameter families of smooth maps of the interval
  with homoclinic orbits or heteroclinic cycles (hypotheses H1-6),
  explosions occur in if and only if there is an isolated crossing
  orbit.
\end{theorem}

\begin{theorem}[General explosion classification]\label{T:general}
  Explosions within generic one-parameter families of smooth maps of
  the interval (hypotheses H1-3) are the result of either
  a tangency between stable and unstable manifolds of fixed
  or periodic points or a saddle node bifurcation of a fixed or
  periodic point. 
\end{theorem}

The paper proceeds as follows: In Section~\ref{S:def}, we give basic
definitions of explosions and homoclinic tangencies. We are
considering the particular recurrence class of chain recurrent points,
which are defined in this section as well. In Section~\ref{S:homo}, we
prove the Explosions at tangencies theorem.  In
Section~\ref{S:mainproof}, we prove the General explosion
classification theorem.  Our results rely on a very sophisticated and
well-developed theory for the dynamics of interval maps.  We have
briefly stated the necessary results in the course of the proof.

\section{Basic definitions\label{S:def}}

We now give some formal definitions of concepts described in the
introduction.  Let $f:(I \times J) \subset (R\times R) \rightarrow R$
be a smooth one-parameter family of maps.  We exchangeably write two
notations: $f(x,\lambda)=f_{\lambda}(x)$.  For the definition of an
explosion it is more natural to use the concept of chain recurrence
rather than recurrence. The relationship between chain recurrence and
other types of recurrence is discussed in~\cite{bonatti:diaz:viana:05}.

\begin{definition} 
For an iterated function $g$, there is an {\bf $\epsilon$-chain} from
$x$ to $y$ when there is a finite sequence $\left( z_0,z_1,\dots,z_N \right)$
such that $z_0=x$, $z_N=y$, and $d(g(z_{n-1}),z_n)<\epsilon$ for all $n$.

If there is an $\epsilon$-chain from $x$ to itself for every
$\epsilon>0$ (where $N>0$), then $x$ is said to be {\bf chain
recurrent}~\cite{bowen:70a,conley:78a}. The {\bf chain recurrent set}
is the set of all chain recurrent points. 
For a one-parameter family $f_\lambda$, we say $(x,\lambda)$ is chain
recurrent if $x$ is chain recurrent for $f_\lambda$.

If for every $\epsilon>0$, there is an $\epsilon$-chain from $x$ to
$y$ and an $\epsilon$-chain from $y$ to $x$, then $x$ and $y$ are said to
be in the same {\bf chain component} of the chain recurrent set. 
\end{definition}

The chain recurrent set and the chain components are invariant under
forward iteration.  The following definition is for an explosion
bifurcation in the chain recurrent set. Such a definition can be
formulated for the non-wandering set as well. 

\begin{definition}[Chain explosions] 
A {\bf chain explosion point} $(x,\lambda_0)$ is a point such that $x$
is chain recurrent for $f_{\lambda_0}$, but there is a neighborhood
$N$ of $x$ such that on one side of $\lambda_0$ (i.e.  either for all
$\lambda<\lambda_0$ or for all $\lambda>\lambda_0$), no point in $N$
is chain recurrent for $f_\lambda$. (All explosion points in this paper 
are chain explosions, so we sometimes drop the qualifier chain.)
\end{definition}

\begin{remark}
  In order the show that $(y,\lambda)$ is {\em not} chain recurrent,
  it is sufficient to show that $y$ is in the closure of the
  hyperbolic periodic orbits for $f_{\lambda}$.
\end{remark}

\begin{remark} In the above definition, at $f_{\lambda_0}$, $x$ is not
necessarily an isolated point of the chain recurrent set. For example,
at a saddle node bifurcation on an invariant circle, the chain
recurrent set consists of two fixed points prior to bifurcation and
the whole circle at and in many cases after bifurcation.
\end{remark}

The chain recurrent set is not invariant under backwards iteration of
a noninvertible map. Thus explosion points are not preserved under
iteration, forward or backward. The following lemma states what is
guaranteed by the fact that chain recurrence is preserved under
forward iteration.

\begin{lemma}\label{L:preimage}
Let $(x,\lambda_0)$ be a chain explosion point for
  $f$. Specifically, there exists $\delta>0$ such
  that there is no chain recurrent point in $B_{\delta}(x)$ for all
  $\lambda<\lambda_0$, but $x$ is chain recurrent at $\lambda_0$. Then
  $f(x)$ is chain recurrent at $\lambda_0$, but $f(x)$ may also be
  chain recurrent for $\lambda<\lambda_0$. In contrast, if 
  $x_{-1}$ is a preimage of $x$, then there is a $\delta_{-1}>0$ such that no
  point in $B_{\delta_{-1}}(x_{-1})$ is chain recurrent for all
  $\lambda<\lambda_0$. Note that $x_{-1}$ may not be chain recurrent at
  $\lambda_0$.
\end{lemma}

\begin{figure}
\begin{center}
\includegraphics[width=0.3\textwidth]{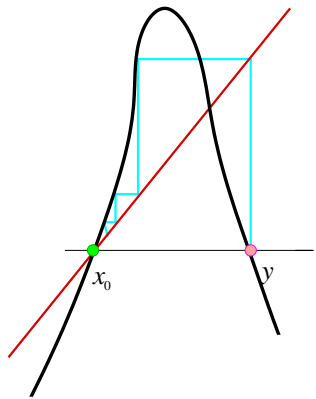}
\caption{A one-dimensional map with a repelling fixed point $x_0$ where
  $f'(x_0)>0$. The point $y$ is a homoclinic point, since
  $f^k(y)=x_0$ (for $k=1$), and there exists a sequence of successive
  preimages of $y$ converging to $x_0$.}
\label{F:fig1}
\end{center}
\end{figure}

We now give definitions of homoclinic and heteroclinic points. Note
that for a diffeomorphism, homoclinic and heteroclinic orbits require
the existence of saddle points with stable and unstable manifolds of
dimension at least one. Thus they can only occur in dimension two or
greater. However, for noninvertible maps, it is possible to have fixed
or periodic points with one-dimensional unstable manifolds and a
zero-dimensional stable manifolds. Marotto terms such points snap-back
repellers~\cite{marotto:78a}. It is not possible to reverse
these stable and unstable manifold dimensions; the existence of a
homoclinic orbit to an attracting fixed point requires a multivalued
map~\cite{sander:00}. The following definition of homoclinic points
for interval maps is depicted in Figure~\ref{F:fig1}.

\begin{definition}[Homoclinic points]
Let $f:R \to R$ be a smooth function with a repelling fixed point
$x_0$. Let $y$ be a point in the unstable manifold of $x_0$ and an
integer $K>0$ such that $f^K(y)=x_0$. Then $y$ is a homoclinic point
to $x_0$. If $x_0$ is periodic with least period $m$, then the same
definition applies by replacing $f$ with $f^m$.
\end{definition}

\begin{remark} If $p$ is a hyperbolic fixed point, and the forward
  limit set of a point $z$ in the unstable manifold of $p$ includes
  $p$, then there are two possibilities: (1) a finite iterate of $z$
  is equal to $p$, (2) the limit set of $z$ contains points other than
  $p$. In this second case, it is not possible for the limit set of
  $z$ to consist of $p$ alone. Therefore $z$ is not referred to as a
  homoclinic point, since its limit set is larger than just $p$. This
  case is considered in later sections of this paper.
\end{remark}

For diffeomorphisms, all orbits through homoclinic points are
homoclinic orbits. For one-dimensional maps, there may be many
non-homoclinic orbits through a homoclinic point.

\begin{definition}[Homoclinic orbits]
  Let $f, x_0,$ and $y$ be as in the above definition of a homoclinic
  point. An orbit $\left( z_{-k} \right)_{k=0}^\infty$ is a
  homoclinic orbit through $y$ if the following conditions are
  satisfied: $z_0=x_0$, $z_{-K}=y$ for some $K$, and for all $k \in
  N$, $f(z_{-k})=z_{-k+1}$, and $\lim _{k \to \infty} z_{-k}=x_0$.
\end{definition}

\begin{figure}
\begin{center}
\includegraphics[width=0.8\textwidth]{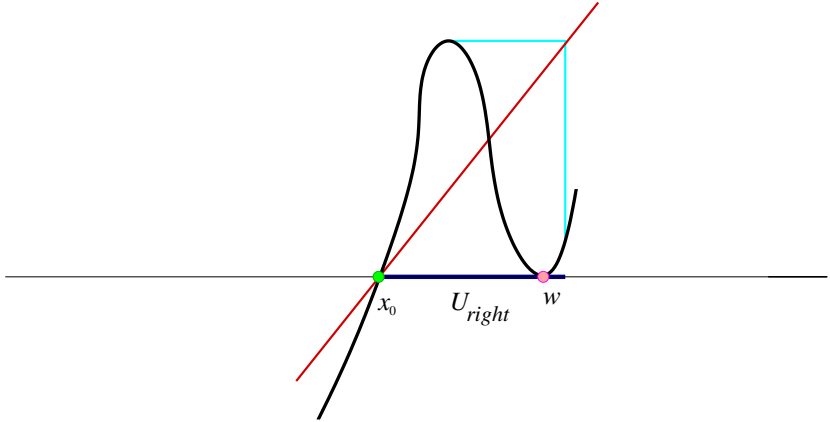}
\caption{A one-dimensional map with a repelling fixed point $x_0$ with
  $f'(x_0)>0$. The point $w$ is a homoclinic tangency point.
  For the particular tangency depicted here, $w$ is contained in a
  {\it non-crossing orbit} (Definition~\ref{D:crossing}), and thus by
  Theorem~\ref{T:NCB} $w$ is not an explosion point.}
\label{F:fig2}
\end{center}
\end{figure}

\begin{figure}
\begin{center}
\includegraphics[width=0.4\textwidth]{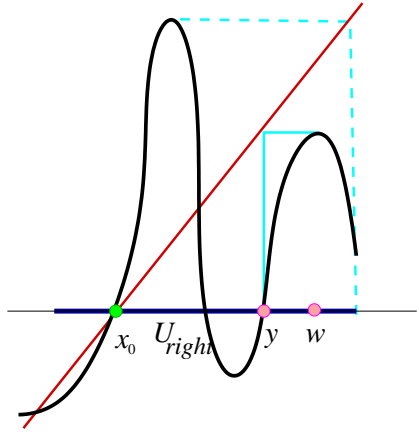}
\caption{In a homoclinic orbit as in Figure~\ref{F:fig2}, the 
tangency point $w$ does not need to be the immediate preimage of the 
fixed point. Here, $w$ is the second preimage of $x_0$.}
\label{F:fig3}
\end{center}
\end{figure}

Since the stable manifold of a homoclinic point is zero-dimensional, a
homoclinic tangency is a tangency of the graph of the map at a
homoclinic point. Homoclinic tangencies are depicted in
Figures~\ref{F:fig2} and~\ref{F:fig3}.

\begin{definition}[Homoclinic tangencies]
Let $f, x_0, y$, and $\left( z_{-k} \right)_{k=0}^\infty$ be as in
the above definition of a homoclinic orbit. The point $w=z_{-L}$ is a
homoclinic tangency point if the graph of $f$ is tangent to the
horizontal line at $w$.
\end{definition}

\section{Explosions at homoclinic tangencies\label{S:homo}}

This section classifies explosions occuring via homoclinic tangencies.
The results are stated for fixed points, but the same results hold for
homoclinic orbits for periodic points with least period $m$ if $f$ is
replaced by $g=f^m$. Throughout this section, we make the following
hypotheses. The first hypothesis is a smoothness assumption.  The
second and third hypotheses are generic assumptions for one-parameter
families. The fourth is a notational convention for the existence of a
homoclinic orbit.  The fifth hypothesis is generic for one-parameter
families containing a homoclinic orbit.

 \begin{enumerate}

 \item[\bf H1] $f:(I \times J) \subset (R \times R) \mapsto I$ is a
   $C^1$ smooth family of $C^2$ interval maps. For a fixed parameter
   $\lambda_0$, denote $f_0:=f_{\lambda_0}$.

 \item[\bf H2] Assume H1. Assume that there are no intervals on which
  $f_0$ is constant.

 \item[\bf H3] Assume H1. Assume that for $f_0$, there is
  at most one of the following:

 \begin{enumerate}
  \item A non-hyperbolic fixed point or periodic orbit.  We assume
    generic behavior as a parameter is varied. That is, non-hyperbolic
    point is a codimension one period doubling or saddle-node
    bifurcation.

  \item One critical point which comprises a tangency between stable
    and unstable manifolds of fixed points or periodic orbits.  We
    assume generic behaviour as a parameter is varied. That is, with
    variation of parameter, the critical point moves from above to
    below (or below to above) the placement of the periodic point.
  \end{enumerate}

 \item[\bf H4] Assume H1, and that for each $\lambda$, $x_{\lambda}$ is a
   repelling fixed point for $f_{\lambda}$. Denote
   $x_0:=x_{\lambda_0}$. Assume that for $f_0$, $y$ is homoclinic to
   $x_0$.

 \item[\bf H5] Assume H1 and H4. At $\lambda=\lambda_0$, the
   homoclinic orbit containing $y$ contains only one critical point.

\end{enumerate}

 The following theorem says that if at a repelling periodic point a
 map has negative derivative then no homoclinic points are explosion
 points.

\begin{theorem}[No explosions for orbits with negative derivatives]
\label{T:negative}
Assume that $f_{\lambda}$ is a family of maps satisfying H1-5, and
$f'_{\lambda_0}(x_0)<0$, then $(y,\lambda_0)$ is not an explosion
point.  
\end{theorem}

\begin{proof}
  We show that
  $y$ is in contained in the closure of the hyperbolic periodic points
  of $f_{\lambda_0}$.

  If $y$ is contained in a homoclinic orbit without a tangency, then
  the homoclinic orbit is preserved under perturbation, which
  automatically implies that $(y,\lambda_0)$ is not an explosion
  point. Thus we assume that $y$ is contained in a homoclinic orbit
  containing a homoclinic tangency point $w$. As before, denote this
  orbit by $\left( z_{-k} \right)_{k=0}^\infty$, where $z_0=x_0$,
  $\lim_{k \to \infty} z_{-k}= x_0$, and $K$ and $L$ are such that
  $y=z_{-K}$ and $w=z_{-L}$.

  Fix a neighborhood $U$ of $y$.  By H5,
  there exists a sequence of neighborhoods $U_{k} \ni z_{-k}$, such
  that: {\em (i)} $U_K \subset U$, {\em (ii)} for $k \ge L$,
  $f_0(U_{k+1})=U_{k}$ and the map is injective, and {\em (iii)} for
  $0 \le k <L$, $U_{L-k}=f_0^k(U_L)$.  Let $V_1=f_0^L(U_L \setminus w)
  = U_0 \setminus x_0$. By H2, $V_1$ is an interval on one side of
  $x_0$. By assumption $f_0'(x_0)<0$, so if the $\left\{ U_k \right\}$
  are sufficiently small, $V_2= f_0^{L+1}(U \setminus w)$ is an
  interval on the other side of $x_0$. That is, $V_1 \cup V_2 \cup
  \left\{x_0\right\}$ is a neighborhood of $x_0$.

  For any sufficiently large $J$, $U_{J} \subset V_i$, where $i \in
  \left\{1,2 \right\}$. In addition, $f^{M}(U_{J})=V_i$, where $M=J$
  or $M=J+1$. Thus for $0 \le k \le J$, $U_{k}$ contains a periodic
  point. Thus there is a periodic point in $U$. Since $U$ was
  arbitrary, there are periodic points $p_j$ with period $n_j$ such
  that $\lim_{j \to \infty} p_j=y$ and $\lim_{j \to \infty} n_j =
  \infty$. By H3, for $j$ sufficiently large, the periodic points are
  hyperbolic. Therefore $(y,\lambda_0)$ is not an explosion point.
\end{proof}

Now consider the positive derivative case.  Let $w$ be a homoclinic
tangency point contained in homoclinic orbit $\left( z_{-k}
\right)_{k=0}^\infty$.  Since the eigenvalue of $x_0$ is positive,
there is a $M$ sufficiently large such that $\left( z_{-k}
\right)_{k=M}^\infty$ lies entirely on one side of $x_0$. That is,
the homoclinic orbit converges to the fixed point along one branch of
the unstable manifold. We denote this by saying that the homoclinic
orbit is contained in the local right or left branch of $x_0$, as
formalized in the following definition.

\begin{definition}[Two branches of the unstable manifold]
  Let $x_0$ be a repelling fixed point for $f \in C^2$, with
  $f'(x_0)>0$. Locally, the left and right branches of the
  unstable manifold of $x_0$ are disjoint. Define $U_{left}$ and $U_{right}$
  to be the respective unions of images of local left and right
  manifold branches.
\end{definition}

\begin{figure}
\begin{center}
\includegraphics[height=0.4\textwidth]{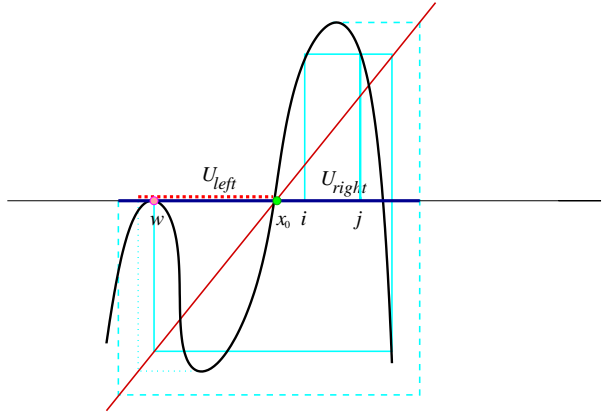}
\caption{The tangency point $w$ is contained in $U_{left} \cap
  U_{right}$, and is thus not an explosion point. However, the
  preimages $i$ and $j$ of $w$ are explosion points.}
\label{F:fig4}
\end{center}
\end{figure}

\begin{figure}
\begin{center}
\includegraphics[height=0.4\textwidth]{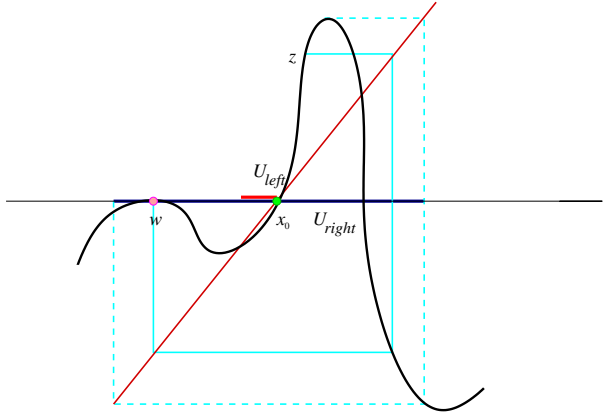}
\caption{In this figure, $w$ is a homoclinic tangency point contained in a
  crossing orbit. Thus $w$ lies in $U_{right}$, but not
  $U_{left}$. By Theorem~\ref{T:crossing}, $w$ is an
  explosion point. As depicted here, the preimage $z$ of $w$ is also
  an explosion point.}
\label{F:fig5}
\end{center}
\end{figure}

\begin{remark} The union of $U_{left}$ and $U_{right}$ is the entire
  unstable manifold of $x_0$.  By the intermediate value theorem
  $U_{left}$ and $U_{right}$ are intervals. If $(U_{left} \cap
  U_{right}) \setminus x_0 $ is not empty, then the intersection must
  contain either $U_{left}$ or $U_{right}$.  For example, $U_{right}$
  may contain points both to the left and to the right of
  $U_{left}$. See Figure~\ref{F:fig4} and~\ref{F:fig5}.
\end{remark}

From the proof of Theorem~\ref{T:negative}, it is clear that to study
chain explosions in homoclinic orbits, it is sufficient to consider
homoclinic orbits containing tangency points. We formalize the
notation in the following hypothesis:

\begin{enumerate}
\item[\bf H6]
For a family satisfying H1 and H4, at $f_0$, point $w$
is a homoclinic tangency point to $x_0$ contained in (at least one)
homoclinic orbit $\left( z_{-k} \right)_{k=0}^\infty$. Let $L$ be
such that $w=z_{-L}$.
\end{enumerate}

\begin{theorem}[No explosions when manifold branches
  intersect]\label{T:branch}
  Assume H1-6. Assume that $f_0'(x_0)>0$ and that $w
  \in U_{right}$. If for any neighborhood $N \ni w$, the 
  sets $f_0^L(N)$ and $U_{right}
  \setminus x_0$ contain points in common, then $(w,\lambda_0)$ is not an
  explosion point.
\end{theorem}

\begin{proof}
  The details of this proof are similar to the negative derivative
  case. Preimages of any small neighborhood $N$ of $y$ are contracting
  and in the local righthand manifold branch. Thus the $L^{th}$ image
  of $N$ includes a shrunk preimage of $N$. Therefore $N$ contains a
  periodic point, which by H3 and H5 is hyperbolic. Since $N$ is
  arbitrary, $(w,\lambda_0)$ is not an explosion point.
\end{proof}

\begin{remark} The theorem above is also true when $U_{right}$ is 
replaced by $U_{left}$. 
\end{remark}

\begin{theorem} \label{T:intersection} Assume H1-6, and that $w$ is
  contained in $U_{right} \cap U_{left}$.  Then $(w,\lambda_0)$ is not
  an explosion point.
\end{theorem}

\begin{proof}
  The image of a neighborhood of $w$ under $f_0^{L}$ contains either
  points of $U_{right}$ or $U_{left}$. But $w$ is contained in both
  $U_{left}$ and $U_{right}$, so by Theorem~\ref{T:branch},
  $(w,\lambda_0)$ is not an explosion point.
\end{proof}

As mentioned in the introduction, in our previous work we gave a
useful geometric method of approaching chain explosions in homoclinic
orbits in two and three dimensions, termed crossing cycles.  In two
dimensions, we showed that a crossing cycle is necessary and
sufficient for a chain explosion to
occur~\cite{alligood:sander:yorke:02}. The analoguous statements are
true in one dimension, as in the theorem below. Crossing and
non-crossing orbits are shown in Figures~\ref{F:fig5} and~\ref{F:fig2}
respectively.

\begin{definition}[Crossing orbits]\label{D:crossing}
  Assume H1-6. Note that for $g:=f_0^{L-1}$, $z_{-1}$ is a point of
  tangency.  If for sufficiently large $k$, {\em (i)} $z_{-k} \in
  U_{left}$, and the graph of $g$ is locally {\em above} the
  horizontal line at $z_{-1}$, or {\em (ii)} $z_{-k} \in U_{right}$,
  and the graph of $g$ is locally {\em below} the horizontal line at
  $z_{-1}$, then we call the homoclinic orbit $\left( z_{-k}\right)$
  a {\bf crossing orbit}. A homoclinic orbit that is not crossing is
  called a {\bf non-crossing orbit}. 
\end{definition}

\begin{theorem}[Explosions imply crossing orbits]\label{T:NCB}
  If H1-6 hold and $(w,\lambda_0)$ is an explosion point, then
  $\left( z_{-k}\right)_{k=0}^{\infty}$ is a crossing orbit.
\end{theorem}

\begin{proof}
  If the graph of $g$ is locally above the horizontal line at
  $z_{-1}$, then $f_0^{L}$ of any neighborhood of $w$ is contained in
  $U_{right}$. Likewise, if the graph is locally below the line,
  $f_0^{L}$ of a neighborhood of $w$ is contained in $U_{left}$.
  The result now follows from Theorem~\ref{T:branch}.
\end{proof}

The results stated so far give necessary conditions for an explosion
point.  We now give a converse to these.  We first make a generic
hypothesis.

By Theorem~\ref{T:intersection}, to classify explosion points, it
suffices to only consider the case of a tangency point $w$ contained
only in one manifold branch of $x_0$, which we state in the following
hypothesis.

\begin{enumerate}
\item[\bf H7] Assume H1-6. Assume that at $\lambda_0$ the tangency
  point $w$ is contained exclusively in one manifold branch of $x_0$. 
\end{enumerate}

\begin{remark}
  Assume H1-7. Then at $\lambda=\lambda_0$, either tangency point $w
  \in U_{right}$, and $U_{left}$ contains no tangencies to fixed points
  or periodic orbits, or tangency point $w \in U_{left}$, and
  $U_{right}$ contains no tangencies to fixed points or periodic
  orbits. This follows from the fact that $w$ is only contained in 
  one of $U_{left}$ and $U_{right}$,
  so by H3 there is no tangency in the other manifold branch.
\end{remark}

The following theorem gives sufficient conditions for an explosion.

\begin{theorem}[Crossing bifurcations and explosions]\label{T:crossing}
  Assume H1-7. If every homoclinic orbit containing $w$ is a crossing
  orbit, then $(w,\lambda_0)$ is an explosion point.
\end{theorem}

\begin{proof}
  For specificity, let $w \in U_{right}$. Since $w$ is only contained in
  crossing orbits, $w$ is not contained in $U_{left}$, and the $L^{th}$
  image of any small neighborhood of $w$ is contained only in $U_{left}$.
  This implies that $U_{left} \subset U_{right}$, but that the subset is
  strict.  As a result, prior to tangency, $w$ is contained in
  $U_{left}$ but not $U_{right}$.  Thus the right endpoint of
  $U_{left}$ is $x_0$. The left endpoint $x_L$ of $U_{left}$ cannot
  map to $x_0$, since that would imply a homoclinic tangency in $U_{left}$.  In
  fact, no points in $U_{left} \setminus x_0$ map to $x_0$, since
  there are assumed to be no tangencies and no intersections with
  $U_{right}$.  If $x_L$ maps into $U_{left}$, then there is a
  neighborhood of $x_L$ which also maps into $U_{left}$.  If $x_L$ is
  a fixed point, then by H3 and H7, it is hyperbolic. It has no
  critical points to mapping to it from the interior of
  $U_{left}$. Thus $x_L$ is an attracting fixed point.  Prior to
  tangency, there is an $\epsilon>0$ such that no $\epsilon$-chains
  carry points from $U_{left}$ to $U_{right}$, and thus from $w$ to
  itself. Thus $(w,\lambda_0)$ is an explosion point.
\end{proof}

We are interested not only in explosion points which are themselves
tangency points, but also in  explosion points far
from tangencies which are caused by tangencies. 
Theorem~\ref{T:negative} contained such a result, but the
other results were specifically about the tangency points. It is now
straightforward to combine the previous results with
Lemma~\ref{L:preimage} to get results for general homoclinic points.
Since the chain recurrent set is invariant under forwards iteration,
the image of non-crossing tangency point is a non-explosion
point. However, there may be explosion points with iterates that are
non-explosion points.  For example, Figure~\ref{F:fig4} shows points of
explosion which are not tangency points but are preimages of tangency
points.

\begin{theorem}[Explosions when manifold branches do not intersect]
  Assume H1-6. If for all $k$, $(z_{-k},\lambda_0)$ is an explosion
  point, then the orbit $\left( z_{-k} \right)_{k=0}^{\infty}$ is a
  crossing orbit. \label{T:noncrossing}
\end{theorem}

\begin{proof}
If the orbit is a non-crossing orbit, then there are hyperbolic periodic 
orbits limiting on every point in the orbit. 
\end{proof}

\begin{theorem}\label{T:preimageinv}
  If H1-6 hold, and $(w,\lambda_0)$ is an explosion point, then any
  pre-image $z_{-k}$ in the homoclinic orbit of $w$ is such that
  $(z_{-k},\lambda_0)$ is an explosion point. 
\end{theorem}

\begin{proof}
  If $(w,\lambda_0)$ is an explosion point, then by
  Lemma~\ref{L:preimage}, $z_{-k}$ is not chain recurrent for
  $\lambda<\lambda_0$. At $\lambda_0$, $z_{-k}$ is contained in a
  homoclinic orbit, and thus chain recurrent. 
\end{proof}

\begin{remark}
This theorem only mentions preimages of $w$ which are contained in a
homoclinic orbit at $\lambda_0$. There may be other preimages of $w$
which are not contained in any homoclinic orbit, and are thus not
chain recurrent at $\lambda_0$.
\end{remark}

\begin{theorem}
  Assume H1-6. Assume that for $\lambda=\lambda_0$, $\left( z_{-k}
  \right)_{k=0}^{\infty}$ is a crossing orbit, but $w$ is also
  contained in some other homoclinic orbit $\left( \alpha_{-k}
  \right)_{k=0}^{\infty}$ which is a non-crossing orbit. Thus
  $(w,\lambda_0)$ is not an explosion point.  
  Then the following statements hold:
\begin{enumerate}
  \item[i.]
  For all $m>0$, $(f^m(w),\lambda_0)$ is not an explosion point. 
  \item[ii.]
  Assume H6 and H7. If there exists $J>L$ such that  $z_{-J}$ is contained neither 
  in $\left( \alpha_{-k} \right)_{k=0}^{\infty}$  nor in any other non-crossing 
  orbit, then for all $k>J$, $(z_{-k},\lambda_0)$ is an explosion point.
\end{enumerate}
\end{theorem}

\begin{proof}

  The first statement holds since at $\lambda_0$, any image of $w$ is contained in
  $\left( \alpha_{-k} \right)_{k=0}^{\infty}$, which is a
  non-crossing orbit. Thus by Theorem~\ref{T:noncrossing}, this image is
  not an explosion point.

  Let $J>L$ be as in Part {\em ii} of the theorem.  For specificity,
  assume that $z_{-J} \in U_{right}$.  Since $z_{-J}$ is contained in
  a crossing orbit, it is not contained in $U_{left}$. A neighborhood
  of $z_{-J}$ maps into $U_{left}$ under $f_0^J$.  If the tangency
  point $w$ is in $U_{right}$, the rest of this proof is exactly the
  same as the proof of Theorem~\ref{T:crossing}, and we conclude that
  $(z_{-J},\lambda_0)$ is an explosion point.  Even if $w \in
  U_{left}$, if the forward orbit of $w$ is in the interior of
  $U_{left}$, the proof follows as before. Furthermore the forward
  orbit of $w$ cannot contain $x_L$, the left endpoint of $U_{left}$:
  Namely, since $x_L$ is the endpoint of an invariant interval, this
  either implies a double homoclinic tangency in the forward orbit of
  $w$, which is ruled out by H3; or it implies both a tangency at $w$
  and a transverse homoclinic intersection at $x_L$ in the orbit of
  $w$, which would mean that $z_{-J}$ is contained in a non-crossing
  orbit.  Finally, using the same proof as given for
  Theorem~\ref{T:preimageinv}, for all $k>J$, $(z_{-k}, \lambda_0)$ is
  an explosion point.
\end{proof}

This completes the classification of when homoclinic points are
explosion points. 

Unlike the planar case, in one dimension the heteroclinic case reduces
to the homoclinic case. That is, if there is a transverse heteroclinic
cycle including an orbit from $\left\{p \right\}$ to $\left\{ q
\right\}$, which are hyperbolic periodic orbits, then both periodic orbits must
be repellers.  Further, the unstable manifold of $\left\{ p \right\}$
contains the unstable manifold of $\left\{ q \right\}$. Since we
assume only one tangency, a heteroclinic tangency point is also a homoclinic
tangency point.

\section{General explosion classification\label{S:mainproof}}

There are many previous results on the structure of $\omega$-limit
sets and chain recurrent sets for interval maps. Block and
Coppel~\cite{block:coppel:92} showed that the chain recurrent set for
maps of the interval can be classified as the set of points $\left\{ x
  : x \in Q(x)\right\}$, where $Q(x)$ is the intersection of all
asymptotically stable sets containing the limit set $\omega(x)$. They
showed that $Q(x)$ is either an asymptotically stable periodic orbit,
a set of asymptotically stable iterated intervals, or special type of
set known as a solenoidal set. However, this classification is not as
useful as it would appear, since $Q(x)$ is not the set of points in
the chain component containing $x$.  Block~\cite{block:78} also proved
that an interval map has a homoclinic point if and only if it has a
periodic point with period not a power of two. Block and
Hart~\cite{block:hart:82} improved on this result to show the
existence of a homoclinic point to a given power of two implies a
cascade of homoclinic bifurcations. Further, if a family of maps
changes from zero to positive entropy, then there is a cascade of
homoclinic bifurcations. Of most relevance to the topic of this paper
are the works of Sarkovskii, Ma\~n\'e~\cite{mane:85},
Blokh~\cite{blokh:98, blokh:95}, and Blokh, Bruckner, Humke, and
Smital~\cite{blokh:etal:96}, where the detailed structure of all
possible $\omega$-limit sets is studied. The $\omega$-limit set can be
a Cantor set known as a basic set (Definition~\ref{D:basic}).  Another
interesting case occurs when the $\omega$-limit set of a point is a
limit of a period doubling cascade, known as a solenoid
(Definition~\ref{D:solenoid}). See~\cite{alseda:etal:98,jimenez:02}
for a detailed characterization of solenoids.  Blokh~\cite{blokh:95}
showed that for $C^2$-smooth maps, $\omega$-limit sets
are either periodic orbits, periodic transitive intervals, subsets of
basic sets, or solenoids.  We use this result to systematically show
that for all possible explosions, there are saddle-node or tangency
bifurcations. 

The following definition describes points at which jumps in an $\epsilon$-chain
are required. We call these points barricades, as they serve to obstruct
orbits.  For example, at the bifurcation parameter,
a saddle node point blocks the points on one side from
reaching the other side.

\begin{definition}[Barricade]\label{D:barricade}
  Assume H1 and H2. Let $z$ be any point. Let
  $S_{\epsilon}=\omega(B_{\epsilon}(z))$. A point $y \in
  \lim_{\epsilon \to 0} S_{\epsilon}$ is called a {\bf barricade} for $z$ if
  it is blocking the orbit of $z$. That is, let $Z_{\epsilon} =
  \omega(B_{\epsilon}(y))$. Then $\lim_{\epsilon \to 0} Z_{\epsilon}$
  contains points not contained in $\lim_{\epsilon \to 0}
  S_{\epsilon}$.
\end{definition}

We consider a point such that the $\omega$-limit set is a fixed point or periodic orbit.

\begin{theorem}\label{T:toper}
  Assume H1-3. Assume that for $f_0$, $z$ is a  point such that
  $\omega(z)$ is equal to a fixed point or periodic orbit $\left\{ p  \right\}$,
  and $\left\{ p  \right\}$  is a barricade.  Then $\left\{ p  \right\}$ 
  is either non-hyperbolic or is the image of a critical point.
\end{theorem}

\begin{proof}
  Assume that $p$ is a fixed point, since otherwise we can let
  $g=f_0^n$. Since $p$ is a barricade, $p \in S_{\epsilon}$, and
  $\lim_{\epsilon \to 0} Z_{\epsilon} \ne \lim_{\epsilon \to 0}
  S_{\epsilon}$ (using the notation from
  Definition~\ref{D:barricade}). This implies that $p$ cannot be an
  attracting fixed point, since for an attracting fixed point
  $\lim_{\epsilon \to 0} Z_{\epsilon}=p$. Thus if $p$ is hyperbolic,
  it must be repelling, meaning that $z$ is a preimage of $p$. Define
  $K$ by $f^K(z)=p$. Since $p$ is a barricade, for small $\epsilon$,
  $f^K(B_{\epsilon}(z))$ is an interval on one side of $p$, implying that $p$ is
  in the orbit of a critical point.
\end{proof}

The previous theorem only indicates that a critical point exists. 
Is still remains to be shown that the critical point is
actually a homoclinic or heteroclinic point.
We use the fact that all $\omega$-limit sets for interval maps
have been classified. 
The following theorem is useful in the proof of several subsequent theorems. 

\begin{theorem}\label{T:useful}
  Assume H1-3, and that $M$ is an invariant interval under $f_0$, such
  that for all $\epsilon>0$, there is an $\epsilon$-chain from a point
  $x \in M$ to a point $y \notin M$. Then there is an $\epsilon$-chain
  from an endpoint $e$ of $M$ to $y$, where $e$ is fixed or period
  two. Furthermore, if $x$ is not an endpoint of $M$, then either $e$
  is non-hyperbolic; or $e$ is repelling, and there is an orbit of a
  critical point in $M$ mapping onto $e$.
\end{theorem}

\begin{proof}
  Since $f(M)=M$, the only way for an $\epsilon$-chain to exit $M$ is
  through an $\epsilon$-jump across one of the endpoints. Call the
  endpoint $e$. Thus there is a chain from $e$ to $y$. Assume $e$ does
  not map to itself or to the other endpoint of $M$. Then $f(e)$ is
  contained in the interior of $M$. But this means that no small
  $\epsilon$-jump at $e$ exits $M$, which is a contradiction. Thus $e$
  can be chosen to be either a fixed point or a period two point.
 
  The orbit of $e$ cannot be attracting, since then there would not be
  $\epsilon$-chains from $e$ to any other point. If the orbit if $e$
  is hyperbolic, then it is repelling, and there is a chain from $x$
  to $e$. If the orbit of $x$ includes $e$, and $x$ is not an endpoint
  of $M$, then the orbit of $x$ contains a critical point, since $M$
  is invariant. If the orbit of $x$ does not include $e$, then the
  limit set of $x$ contains more than just a repelling periodic orbit.
\end{proof}

The following result shows what happens if the backward limit set of a
point contains a periodic orbit.

\begin{theorem}\label{T:fromper}
  Assume H1-3. Let $\left\{ p \right\}$ be a fixed point or
  periodic orbit for $f_0$ which is hyperbolic. Let $z$ be a point
  such that for all $\epsilon>0$ there is an $\epsilon$-chain 
  from $\left\{ p \right\}$ to
  $z$, but $z$ is not in the unstable manifold of $\left\{ p
  \right\}$. Then there is a homoclinic tangency to a periodic orbit.
\end{theorem}

\begin{proof}
  Assume without loss of generality that $p$ is a fixed point (since
  otherwise, we can use the same proof for an iterate of $f_0$). Since
  $p$ is hyperbolic, it is repelling. The unstable manifold of $p$ is
  an invariant interval, denoted $U(p)$. Since $z \notin U(p)$, there
  is a barricade point for $p$. By Theorem~\ref{T:useful}, there must
  be a barricade point which is a periodic endpoint of $U(p)$, with a
  critical point in $U(p)$ mapping to the endpoint. That is, there is
  a homoclinic tangency point for the periodic endpoint.
\end{proof}

Combining Theorems~\ref{T:toper} and~\ref{T:fromper}, we conclude
that if an explosion occurs at a point $(z,\lambda_0)$ such that
$\omega(z)$ is equal to a periodic orbit, then there is either a
saddle-node bifurcation point or a tangency between stable and
unstable manifolds. We now consider more general $\omega$-limit sets. 
First consider the case where the $\omega$-limit set of a point is an interval. 

\begin{theorem}
  Assume H1-3, for $f_0$, $z$ is such that $\omega(z)$ is an
  interval $M$, and that there exists $x \in M$ such that $e$ is a
  barricade for $z$.  Then $e$ is an endpoint of the interval, $e$ is
  fixed or period two, and if $e$ is hyperbolic, then there is a
  homoclinic tangency to $e$.
\end{theorem}

\begin{proof}
  Since $\omega(z)=M$, $f_0$ is transitive on $M$. Therefore the
  unstable manifold of any repelling periodic orbit in $M$ contains
  all of $M$.  A barricade must not be in the interior of
  $\omega(z)$. Thus $e$ is an endpoint of $M$.  The result now follows from
  Theorem~\ref{T:useful}.
\end{proof}

We now consider the case of an $\omega$-limit set that is contained in
an invariant interval but is nowhere dense.

\begin{definition}[Basic set]\label{D:basic} Assume H1.  
  Let $M=\cup_{k=1}^n
  I_k$ be an $n$-periodic cycle of intervals for a function
  $f_0$. Define $B(M,f_0)=\left\{ x \in M: \right.$ for every open
  interval $U$ of $x$ in $M, \overline{\mbox{orb}(U)} = \left. M
  \right\}$. If $B(M,f_0)$ is infinite, then it is called a basic set.
\end{definition}

\begin{theorem}Assume H1, H2, and H3a. Assume that $(z,\lambda_0)$ is an
  explosion point, and $\omega(z)$ is nowhere dense and is contained
  in a basic set $B(M,f_0)$. Then $\omega(z)$ is a periodic
  orbit. 
\end{theorem}

\begin{proof}
  Assume that $(z,\lambda_0)$ is an explosion point, and that
  $\omega(z) \subset B(M,f_0)$. If $z$ is contained in an interval
  complementary to $B(M,f_0)$, then by Blokh~\cite{blokh:98},
  $\omega(z)$ is a periodic orbit.

  Assume $\omega(z)$ is non-periodic, meaning $z$ is contained
  in the basic set. By~\cite{blokh:95}, $B(M,f_0)$ is contained in the
  closure of the periodic orbits for $f_0$. Using H3a, there is a sequence
  of hyperbolic periodic points converging to $z$. Thus $(z,\lambda_0)$
  is not an explosion point.
\end{proof}

By the above theorem combined with Theorems~\ref{T:toper}
and~\ref{T:fromper}, if an explosion point has an $\omega$-limit set
which is a basic set, then there is either a saddle-node point or a
tangency.  The last possibility for an $\omega$-limit set is a
solenoid, as in the following definition.

\begin{definition}[Solenoid]\label{D:solenoid}
  Assume H1. Let $M_j = \cup_{k=1}^{n_j} I_k^j $ be
  a nested sequence of cycles of intervals for a function $f_0$ with
  least period $n_j$, $\lim_{j \to \infty} n_j = \infty$.  Thus
  $f_0^{n_j}(I_1^j)=I_1^j$, $n_j$ is increasing, and for each $j$,
  $I_1^{j+1} \subset I_1^{j}$.  If the set $S=\cap_{j=1}^{\infty} M_j$
  is nowhere dense, then $S$ is called a solenoid or Feigenbaum-like
  set.
\end{definition}

Jim{\'e}nez L{\'o}pez has shown that solenoids are the boundary of
chaos and order~\cite{jimenez:02}.  Blokh~\cite{blokh:95} demonstrated
that solenoids and basic sets are disjoint.  We prove the
following result.

\begin{theorem}Assume H1, H2, and H3a.  Assume that $(z,\lambda_0)$ is an
  explosion point, and $\omega(z)$ is a solenoid $S$.  Then there
  is an infinite sequence periodic orbits which are barricades with 
  associated tangencies.
\end{theorem}

\begin{proof}
  Since solenoids and basic sets are closed, invariant, and disjoint,
  there is a neighborhood of solenoid $S$ containing no basic sets.
  By definition, there is an infinite nested sequence of invariant
  cycles of intervals $M_j$ containing $S$. Blokh~\cite{blokh:98}
  proved that the periodic orbits are dense in a neighborhood of $S$,
  meaning that $z$ is not contained in $S$. Thus for all $j$
  sufficiently large, for all $\epsilon>0$ 
  there is an $\epsilon$-chain from points in
  $M_j$ to $z$, but $z$ is not contained in $M_j$. By
  Theorem~\ref{T:useful}, there exists $e_j$, an endpoint of an
  interval the cycle of $M_j$ which is periodic. By hypothesis H3a, for
  sufficiently large $j$, $e_j$ is hyperbolic. Since for all $\epsilon>0$ 
  there is an
  $\epsilon$-chain from $e_j$ to $z$, $e_j$ is a repeller for large
  $j$. There is an orbit of a critical point in $M_j$ mapping onto
  $e_j$. Furthermore, $z$ is not contained in the unstable manifold of
  $e_j$, since there is a nested sequence of invariant $M_j$ not
  containing $z$. Theorem~\ref{T:fromper} implies that the critical
  point to $e_j$ is a point of homoclinic tangency. 
\end{proof}

\begin{remark}
  If $f_0$ has a finite number of homoclinic and heteroclinic
  tangencies, as assumed in H3b, then the above theorem shows that there are
  no forward chains from solenoid $S$ to a point outside of $S$.
\end{remark}

\begin{remark}\label{R:sn} We have shown that there is either a tangency 
  or a non-hyperbolic critical point contained in the same chain component as
  $z$. Under a generic hypothesis (H3a), a
  non-hyperbolic periodic orbit is either codimension-one saddle-node
  or period doubling bifurcation. In fact, such an orbit is not a
  period doubling point, since the periodic orbit at a period-doubling
  bifurcation point is attracting.
\end{remark}

We now combine the results of this section to give a proof of the
General explosion classification theorem.

\begin{proof} (of the General explosion classification theorem) 
  Assume that $(z,\lambda_0)$ is an explosion point.  The
  only possibilities for $\omega(z)$ are a periodic orbit, a cycle of
  intervals, a nowhere dense basic set, and a solenoid. Above, we have
  shown that in any of these cases, there is a periodic barricade
  point for $z$ which is either non-hyperbolic, or there is a
  homoclinic or heteroclinic tangency. In fact, the case of a solenoid
  is ruled out by H3b. By Remark~\ref{R:sn}, a non-hyperbolic periodic
  orbit must necessarily be a saddle-node point. This completes the
  proof the theorem.
\end{proof}

\section{Acknowledgements}

E.S. would like to thank Alexander Blokh and Victor Jim\'enez L\'opez
for their helpful suggestions. J.A.Y.  was partially supported by NSF
Grant DMS0616585.

\begin{flushleft}
%
%
\addcontentsline{toc}{subsection}{References}
\footnotesize
%
%

\bibliographystyle{abbrv}

\begin{thebibliography}{10}

\bibitem{alligood:sander:yorke:02}
K.~Alligood, E.~Sander, and J.~Yorke.
\newblock Explosions: Global bifurcations at heteroclinic tangencies.
\newblock {\em Ergodic Theory and Dynamical Systems}, 22(4):953--972, 2002.

\bibitem{alligood:sander:yorke:06}
K.~Alligood, E.~Sander, and J.~Yorke.
\newblock Three-dimensional crisis: Crossing bifurcations and unstable
  dimension variability.
\newblock {\em Physical Review Letters}, 96(244103), 2006.

\bibitem{alseda:etal:98}
L.~Alsed{\`a}, V.~Jim{\'e}nez~L{\'o}pez, and L.~Snoha.
\newblock All solenoids of piecewise smooth maps are period doubling.
\newblock {\em Fund. Math.}, 157(2-3):121--138, 1998.
\newblock Dedicated to the memory of Wies\l aw Szlenk.

\bibitem{block:78}
L.~Block.
\newblock Homoclinic points of mappings of the interval.
\newblock {\em Proc. Amer. Math. Soc.}, 72(3):576--580, 1978.

\bibitem{block:coppel:92}
L.~Block and W.~Coppel.
\newblock {\em Dynamics in One Dimension}, volume 1513 of {\em Lecture Notes in
  Mathematics}.
\newblock Springer, Berlin, 1992.

\bibitem{block:hart:82}
L.~Block and D.~Hart.
\newblock The bifurcation of homoclinic orbits of maps of the interval.
\newblock {\em Ergodic Theory Dynamical Systems}, 2(2):131--138 (1983), 1982.

\bibitem{blokh:98}
A.~Blokh.
\newblock Density of periodic orbits in {$\omega$}-limit sets with the
  {H}ausdorff metric.
\newblock {\em Real Anal. Exchange}, 24(2):503--521, 1998/99.

\bibitem{blokh:etal:96}
A.~Blokh, A.~Bruckner, P.~Humke, and J.~Sm\'{i}tal.
\newblock The space of $\omega$-limit sets of a continuous map of the interval.
\newblock {\em Transactions of the American Mathematical Society},
  348(4):1357--1372, 1996.

\bibitem{blokh:95}
A.~M. Blokh.
\newblock The ``spectral'' decomposition for one-dimensional maps.
\newblock In {\em Dynamics reported}, volume~4 of {\em Dynam. Report.
  Expositions Dynam. Systems (N.S.)}, pages 1--59. Springer, Berlin, 1995.

\bibitem{bonatti:diaz:viana:05}
C.~Bonatti, L.~Diaz, and M.~Viana.
\newblock {\em Dynamics beyond hyperbolicity}.
\newblock Springer Verlag, Berlin, 2005.

\bibitem{bowen:70a}
R.~Bowen.
\newblock Markov partitions for axiom a diffeomorphisms.
\newblock {\em American Journal of Mathematics}, 92:725--747, 1970.

\bibitem{conley:78a}
C.~Conley.
\newblock {\em Isolated Invariant Sets and the Morse Index}.
\newblock American Mathematical Society, Providence, R.I., 1978.

\bibitem{diaz:rocha:02}
L.~D{\'{\i}}az and J.~Rocha.
\newblock Heterodimensional cycles, partial hyperbolicity and limit dynamics.
\newblock {\em Fundamenta Mathematicae}, 174(2):127--186, 2002.

\bibitem{horita:muniz:sabini:07}
V.~Horita, N.~Muiz, and P.~R. Sabini.
\newblock Non-periodic bifurcations of one-dimensional maps.
\newblock {\em Ergodic Theory and Dynamical Systems}, Online, 2007.

\bibitem{jimenez:02}
V.~Jim{\'e}nez~L{\'o}pez.
\newblock Period doubling is the boundary of chaos and of order in the {$C\sp
  1$}-topology of interval maps.
\newblock {\em Nonlinearity}, 15(3):817--839, 2002.

\bibitem{mane:85}
R.~Ma{\~n}{\'e}.
\newblock Hyperbolicity, sinks and measure in one-dimensional dynamics.
\newblock {\em Comm. Math. Phys.}, 100(4):495--524, 1985.

\bibitem{marotto:78a}
F.~Marotto.
\newblock Snap-back repellers imply chaos in ${R}^n$.
\newblock {\em Journal of Mathematical Analysis and Applications}, 63:199--223,
  1978.

\bibitem{newhouse:palis:76}
S.~Newhouse and J.~Palis.
\newblock Cycles and bifurcation theory.
\newblock {\em Ast\'{e}risque}, 31(44--140), 1976.

\bibitem{palis:takens:93a}
J.~Palis and F.~Takens.
\newblock {\em Hyperbolicity \& Sensitive Chaotic Dynamics at Homoclinic
  Bifurcations}.
\newblock Cambridge University Press, Cambridge, 1993.

\bibitem{sabini:01}
P.~R. Sabini.
\newblock {\em Non-periodic bifurcations at the boundary of hyperbolic
  systems}.
\newblock PhD thesis, IMPA, 2001.

\bibitem{sander:00}
E.~Sander.
\newblock Homoclinic tangles for noninvertible maps.
\newblock {\em Nonlinear Analysis}, 41(1-2):259--276, 2000.

\end{thebibliography}

{\bf AMS Subject Classification: 37.}\\[2ex]

%
E.~Sander
Department of Mathematical Sciences,
George Mason University,
4400 University Dr.,
Fairfax, VA, 22030, USA.\\
E-mail: \texttt{sander@math.gmu.edu}\\[2ex]
J.A.~Yorke,
IPST,
University of Maryland,
College Park, MD 20742, USA. \\
E-mail: \texttt{yorke@ipst.umd.edu}\\[2ex]

\end{flushleft}

\end{document}